\documentclass[12pt,amsfonts]{article}
\usepackage[margin=1in]{geometry}  
\usepackage{graphicx}              
\usepackage[latin1]{inputenc}
\usepackage{amsmath, amssymb, amsthm, epsfig}               
\usepackage{enumerate}
\usepackage{amsfonts}              
\usepackage{amsthm}                
\usepackage{amsthm}
\usepackage{enumerate}
\theoremstyle{plain}
\usepackage{color}

\usepackage{tikz}  
\usepackage{pgfplots}

\def\nd{\noindent}
\def\thend{\rule{3mm}{3mm}}

\newtheorem{theorem}{Theorem}[section]

\newtheorem{proposition}{Proposition}[section]

\newtheorem{remark}[theorem]{Remark}

\newtheorem*{theorem*}{Theorem}

\numberwithin{equation}{section}

\begin{document}
\title{Nodal solutions for the double phase problems}

\author{ Chao Ji,\footnote{Corresponding author} \footnote{C. Ji was partially supported by National Natural Science Foundation of China (No. 12171152).} \,\, Nikolaos S. Papageorgiou }

\maketitle

\begin{abstract}
We consider a parametric nonautonomous $(p, q)$-equation with unbalanced growth as follows
 \begin{align*}
 \left\{
\begin{aligned}
&-\Delta_p^\alpha u(z)-\Delta_q u(z)=\lambda \vert u(z)\vert^{\tau-2}u(z)+f(z, u(z)), \quad
\quad
\hbox{in }\Omega,\\
&u|_{\partial \Omega}=0,
\end{aligned}
\right.
\end{align*}
where $\Omega \subseteq \mathbb{R}^N$ be a bounded domain with Lispchitz boundary $\partial\Omega$, $\alpha \in L^{\infty}(\Omega)\backslash \{0\}$, $a(z)\geq 0$ for a.e. $z \in \Omega$, $ 1<\tau< q<p<N$ and $\lambda>0$. In the reaction there is a parametric concave term and a perturbation $f(z, x)$. Under the minimal conditions on $f(z, 0)$, which essentially restrict its growth near zero, by employing variational tools,  truncation and comparison techniques, as well as critical groups, we prove that for all small values of the parameter $\lambda>0$, the problem has at least three nontrivial bounded solutions (positive, negative, nodal), which are ordered and asymptotically vanish as $\lambda \rightarrow 0^{+}$.

\end{abstract}

{\small \textbf{2010 Mathematics Subject Classification:} 35J60, 58E05.}

{\small \textbf{Keywords:} Unbalanced growth, Generalized Orlicz spaces, Critical groups, Extremal constant sign solutions, Nodal solutions.}

\section{Introduction}

In this paper we are concerned with the following parametric double phase Dirichlet problem
\begin{align}
\left\{
\begin{aligned}
&-\Delta_p^\alpha u(z)-\Delta_q u(z)=\lambda \vert u(z)\vert^{\tau-2}u(z)+f(z, u(z)), \quad \quad \text{in }\Omega,\\
&u|_{\partial \Omega}=0,
\end{aligned}
\right.
\tag{$P_{\lambda}$}
\end{align}
where $\Omega \subseteq \mathbb{R}^N$ be a bounded domain with Lispchitz boundary $\partial\Omega$, $\alpha \in L^{\infty}(\Omega)\backslash \{0\}$, $a(z)\geq 0$ for a.e. $z \in \Omega$, $ 1<\tau< q<p<N$ and $\lambda>0$. We denote the weighted $r$-Laplacian differential operator by $\Delta_r^a$ and define it as follows
$$
\Delta_r^\alpha u=\operatorname{div}\left(a(z)\nabla u|^{r-2} \nabla u\right) .
$$
where $1<r<\infty$.  If $\alpha \equiv 1$, then we recover the usual $r$-Laplacian differential operator. Problem $(P_{\lambda})$ is driven by the sum of two such operators with different exponents, making it a non-homogeneous differential operator. This operator is related to the double phase energy functional, which is defined by
$$
u \rightarrow \int_{\Omega}\Big(\alpha(z)\vert\nabla u\vert^p+\vert\nabla u\vert^q\Big)dz.
$$
It is noteworthy that we do not assume that the weight function $\alpha(\cdot)$ is bounded away from zero, i.e., we do not require that $\underset{\Omega}{\text{ess}\inf}\,\alpha>0$.
As a result, the density function of the above integral functional, denoted as the integrand $\eta(z, t)=a(z) t^p+t^q$, exhibits unbalanced growth, which can be characterized by:
$$
t^q\leq \eta(z, t)\leq c_{0}(1+t^p),\,\,\text {for a.e.}\,\, z\in\Omega, \,\,\text {all}\,\, t\geq 0,\,\,\text {some}\,\, c_{0}>0.
$$
Such integral functionals were first considered by Marcellini \cite{Ma, Ma1} and Zhikov \cite{Zh1, Zh2} in the context of problems in the calculus of variations (including the Lavrentiev gap phenomenon) and of nonlinear elasticity theory.  For problems with unbalanced growth, only local regularity results exist (please see the survey papers \cite{Ma2} due to Marcellini and \cite{MR} due to Mingione and R\v{a}dulescu),  and there is no global regularity theory (i.e., regularity up to the boundary). This limitation restricts the use of many powerful technical tools available for problems with balanced growth.

In the present paper, we address these challenges and demonstrate, under minimal conditions on the perturbation $f(z, x)$, which essentially impose conditions on $f(z, \cdot)$ only near zero (local perturbation), the nodal solution and multiple solutions of problem $(P_{\lambda})$ are studied.

The main result of this paper is the following.
\begin{theorem}\label{mainth}
If hypotheses $(H_0), (H_1)$ hold (see Section \ref{Sec2}), then for all $\lambda>0$ small, problem $(P_{\lambda})$ possesses at least three nontrivial, bounded solutions with sign information (positive, negative, and nodal/sign-changing). These solutions are ordered and converge to zero in $L^{\infty}(\Omega)$ as $\lambda \rightarrow 0^{+}$.
\end{theorem}

The absence of a global regularity theory makes it challenging and difficult to produce nodal solutions for the double phase problems, some authors have already done some interesting work in this direction. In particular, in \cite{CBW}, Crespo-Blanco and Winkert  considered quasilinear elliptic equations driven by the variable exponent double phase operator with superlinear right-hand sides. Under very general assumptions on the nonlinearity, they proves a multiplicity result for such problems and showed the existence of a positive solution, a negative one and a solution with changing sign; By using variational methods,  Liu and Dai \cite{LD} obtained various existence and multiplicity results for the following double phase problem
$$
\begin{cases}-\operatorname{div}\left(|\nabla u|^{p-2} \nabla u+a(x)|\nabla u|^{q-2} \nabla u\right)=f(x, u) & \text { in } \Omega, \\ u=0 & \text { on } \partial \Omega .\end{cases}
$$
In particular, they found a sign-changing ground state solution; then, Papageorgiou and Zhang \cite{PZ} dealt with a nonlinear unbalanced double-phase problem with a superlinear reaction and Robin boundary condition.  They showed that the problem has three nontrivial solutions all with sign information (positive, negative and nodal). We notice that in order to study nodal solutions, they all applied the Nehari manifold method and its variants in the papers mentioned above.  However, these methods required restrictive monotonicity conditions (see hypothesis $\left(f_6\right)$ of \cite{CBW}, hypothesis $\left(f_5\right)$ of \cite{LD}, hypothesis $H_1($ iii) of \cite{PZ}) or differentiability conditions on $f(z, \cdot)$ (see hypothesis $H_1$  differentiability condition in \cite{PZ}), which we avoid in our study. Instead, in this paper we will employ variational tools, along with truncation and comparison techniques and critical groups, to establish the existence of nodal solutions under minimal hypotheses.

Recently, in \cite{PVV}, Papageorgiou, Vetro and Vetro  considered  problem $(P_{\lambda})$. They proved that for all parametric values $\lambda>\lambda^*$ the problem has at least three nontrivial solutions, two of which have constant sign, one is sign-changing solution, here the critical parameter $\lambda^*$ is precisely in terms of the spectrum of the $q$-Laplacian. Notice that in \cite{PVV}, the perturbation enters with a negative sign and $\tau=q$. In this paper, the perturbation enters with positive sign and $\tau<q$. There is a perturbation of the $q$-eigenvalue problem in \cite{PVV} while in our case we have a problem with a concave term and a general perturbation. That is why, in \cite{PVV}, the focus is on the existence of solutions for large values of $\lambda$, while in this paper, we will instead consider the existence of solutions for small values of $\lambda$. Finally we would like to mention that in contest to \cite{PVV}, here we do not have any asymptotic condition on $\frac{f(z, x)}{\vert x\vert^{p-2}x}$ as $x\rightarrow \pm \infty$. To the best of our knowledge, this is the first work producing nodal solutions for unbalanced growth problems using critical groups. This is rather surprising given the lack of a global regularity theory for such problems.

The paper is organized as follows. In Section \ref{Sec2} we introduce the functional setting and give some preliminaries. Then, in Section \ref{Sec3}, we prove existence of constant sign solutions.  In  Section \ref{Sec4}, we produce the nodal solutions. Finally, in the last section, we will give the proof of Theorem \ref{mainth}.

\section{The variational framework and some preliminaries}\label{Sec2}
As a consequence of the unbalanced growth of the function $\eta(z, \cdot)$, the standard setting of the classical Lebesgue and Sobolev spaces is inadequate. Instead, we need to work with the more suitable generalized Orlicz spaces. For a comprehensive presentation of the theory of these spaces, we refer to the book by Harjulehto and H\"asto \cite{HH}.

Our hypotheses on the weight $\alpha(\cdot)$ and the exponents $p, q, \tau$ are as follows:
\begin{equation}
\label{V}\tag{$H_{0}$}
\alpha \in L^{\infty}(\Omega) \backslash\{0\}, a(z) \geqslant 0 \,\,\text{for a.e.} \,\,z\in\Omega,\, 1<\tau<q<p<N, \alpha\leq q^*=\frac{N q}{N-q}\,\, \text{and} \,\,\frac{p}{q}<1+\frac{1}{N}.
\end{equation}

\begin{remark}\label{regula}
The last inequality implies that the exponents $p, q$ can not be too far apart. Also,  it leads to that $p<q^*$ and this in turn leads to the compact embeddings of some relevant spaces (see Proposition \ref{embeddings} below). The hypothesis on the weight function $\alpha(\cdot)$, together with \cite[Proposition 2.18]{CGHW}, guarantees the validity of the Poincar\'{e} inequality on the generalized Sobolev-Orlicz space $W_0^{1, \eta}(\Omega)$, which we will introduce later.
\end{remark}

Let $L^{0}(\Omega)$ be the space of all measurable functions $u: \Omega \rightarrow \mathbb{R}$. As usual, we identify two such functions which differ only on a Lebesgue null subset of $\Omega$. The generalized Lebesgue-Orlicz space $L^\eta(\Omega)$ is defined by
$$
L^\eta(\Omega)=\left\{u \in L^0(\Omega): \rho_\eta(u)=\int_{\Omega} \eta(z,|u|) d z<\infty\right\},
$$
where the function $\rho_\eta(\cdot)$ is known as the modular function. This space is equipped with the so-called Luxemburg norm $\|\cdot\| \eta$, defined by
$$
\|u\|_\eta=\inf\Big\{ \lambda>0: \rho_\eta\left(\frac{u}{\lambda}\right) \leqslant 1\Big\}.
$$
With this norm, the space $L^\eta(\Omega)$ becomes a separable, reflexive Banach space that is also uniformly convex due to the  uniformly convexity of $\eta\left(z, \cdot\right)$.

Using $L^{\eta}(\Omega)$, we can define the corresponding generalized Sobolev-Orlicz space $W^{1, \eta}(\Omega)$ as
$$
W^{1, \eta}(\Omega)=\Big\{u \in L^\eta(\Omega):|\nabla u| \in L^\eta(\Omega)\Big\},
$$
where $\nabla u$ represents the weak gradient of $u$. We equip $W^{1, \eta}(\Omega)$
with the norm $\|\cdot\|_{1, \eta}$ defined by
$$
\|u\|_{1, \eta}=\|u\|_\eta+\|\nabla u\|_\eta,\,\, \text {for any}\,\,u\in W^{1, \eta}(\Omega),
$$
here $\Vert \nabla u\Vert_\eta=\Vert \vert \nabla u \vert \Vert_\eta$. Additionally,  we define
$$
W_0^{1, \eta}(\Omega)=\overline{C_c^{\infty}(\Omega)}^{\|\cdot\|_{1, \eta}}.
$$
Due to the Poincar\'{e} inequality being valid on $W_0^{1, \eta}(\Omega)$, we can consider the equivalent norm on $W_0^{1, \eta}(\Omega)$
$$
\|u\|=\|\nabla u\|_\eta,  \,\,\text{for all}\,\, u \in W_0^{1, \eta}(\Omega).
$$
Both $W^{1, \eta}(\Omega)$ and $W_0^{1, \eta}(\Omega)$ are separable, reflexive Banach spaces, with the uniformly convexity.

We have the following  embedded results among these spaces  which are useful.

\begin{proposition}\label{embeddings}
\begin{enumerate}
	\item $L^\eta(\Omega)\hookrightarrow L^s(\Omega), W_0^{1, \eta}(\Omega) \hookrightarrow W_{0}^{1, s}(\Omega)$ continuously for all $s\in [1, q]$;
	\item $L^{p}(\Omega) \hookrightarrow L^\eta(\Omega)$ continuously;
  \item  $W_0^{1, \eta}(\Omega) \hookrightarrow L^s(\Omega)$ continuously if $s \in [1, q^*]$, and $W_0^{1, \eta}(\Omega) \hookrightarrow L^s(\Omega)$ compactly if $s \in [1, q^*)$, where $q^* = \frac{Nq}{N - \eta q}$ is the critical Sobolev exponent.
\end{enumerate}
\end{proposition}

Note that the modular function $\rho_\eta: W_0^{1, \eta}(\Omega)\rightarrow\mathbb{R}^{+}$ is continuous and convex, hence by Mazur's lemma, it is weakly lower semi-continuous. There is a close relation between the modular function $\rho_\eta(\cdot)$ and the norm $\Vert \cdot \Vert$ as follows.
\begin{proposition}\label{normre}

\begin{enumerate}
	\item \text{If} $u\neq 0$\text{, then} $\|u\|=\lambda \Leftrightarrow \rho_\eta\left(\frac{\nabla u}{\lambda}\right)=1$;
	\item $\|u\|<1$ (\text{respectively}\, $=1,>1) \Leftrightarrow \rho(\nabla u)<1$(\text{respectively}\, $=1,>1)$;
  \item  $\|u\|<1 \Rightarrow\|u\|^p \leq \rho(\nabla u) \leq \|u\|^q$;
 \item $\|u\|>1 \Rightarrow\|u\|^q \leq \rho(\nabla u) \leq \|u\|^p$;
\item $\|u\| \rightarrow \infty ( \text{respectively} \rightarrow 0) \Leftrightarrow \rho(\nabla u) \rightarrow \infty ( \text{respectively} \rightarrow 0)$.
\end{enumerate}

\end{proposition}

Furthermore, we introduce the map $V: W_0^{1, \eta}(\Omega) \rightarrow (W_0^{1, \eta}(\Omega))^*$ defined by
$$
\langle V(u), h\rangle=\int_{\Omega}\Big(a(z)|\nabla u|^{p-2}+|\nabla u|^{q-2}\Big)(\nabla u, \nabla h)dz, \,\,\text{for any}\,\, u, h \in W_0^{1, \eta}(\Omega).
$$
This map has the following important properties.
\begin{proposition}\label{mappro}
 The map $V(\cdot)$ is bounded, continuous, strictly monotone (thus maximal monotone too) and of type $(S)_{+}$, which means that
$$
 \text {If}\,\, u_n \rightharpoonup u\,\, \text{in}\,\, W_0^{1, \eta}(\Omega)\, \text { and }\, \limsup _{n \rightarrow \infty}\left\langle V\left(u_n\right), u_n-u\right\rangle\leq 0,
$$
then $u_n \rightarrow u$  in $W_0^{1,\eta}(\Omega)$.
\end{proposition}
If $u \in L^0(\Omega)$, we define
$$
u^{ \pm}(z)=\max \{ \pm u(z), 0\},\,\, \text { for a.e.}\,\, z \in \Omega.
$$
Observe that $u=u^{+}-u^{-}$ and $\vert u\vert=u^{+}+u^{-}$. Additionally, if $u \in W_0^{1, \eta}(\Omega)$, then $u^{\pm} \in W_0^{1, \eta}(\Omega)$. Given $h_{1}, h_2 \in L^0(\Omega)$ with $h_1(z) \leqslant h_2(z)$ for a.e. $z \in \Omega$, we define the order interval $[h_1, h_2]$ as
$$
\left[h_1, h_2\right]=\left\{u \in W_0^{1, \eta}(\Omega): h_1(z) \leqslant u(z) \leqslant h_2(z), \,\,\text {for a.e.} \,\,z \in \Omega\right\}.
$$
If $X$ is a Banach space and $\varphi \in C^1(X)$, then
$$
\left.K_{\varphi}=\left\{u \in X: \varphi^{\prime}(u)=0\right\} \text { (critical points set of } \varphi\right) \text {. }
$$

A set $C \subseteq W_0^{1, \eta}(\Omega)$ is said to be "downward directed" (respectively, "upward directed"), if given $u_1, u_2 \in C$, we can find $u\in C$ such that $u \leq u_1, u \leq u_2$ (respectively, if given $v_1, v_2 \in C$, we can find $v \in C$ such that $\left.v_1 \leq v, v_2 \leq v\right)$.

As we already mentioned in the introduction, in order to overcome the serious difficulties arising from the absence of a global regularity theory, we will use the critical groups and their properties. So, let us briefly recall some basic definitions and facts from that theory. For the details we refer to Chapter 6 in the book \cite{PRR}.

Let $X$ be a Banach space and $\left(Y_{1}, Y_2\right)$ be a topological pair where $Y_2 \subseteq Y_1 \subseteq X$. For this pair, $H_k\left(Y_1, Y_2\right)$, $k \in \mathbb{N}_{0}$, denotes the $k$th-relative singular homology group with integer coefficients.  Given $\varphi \in C^{1}(X)$ and $c \in \mathbb{R}$, we define $\varphi^c=\{u \in X: \varphi(u) \leq c\}$. If $u \in K_{\varphi}$ is isolated and $c=\varphi(u)$, then the critical groups of $\varphi(\cdot)$ at $u$ are given by
$$
C_k(\varphi, u)=H_k\left(\varphi^c \cap \mathcal{U}, \varphi^c \cap \mathcal{U}\backslash\{u\}\right),\,\, k \in \mathbb{N}_0,
$$
with $\mathcal{U}$ a neighborhood of $u$ such that $K_{\varphi} \cap \varphi^c \cap \mathcal{U}=\{u\}$. These critical groups are well-defined and independent of the choice of the isolating neighborhood $\mathcal{U}$, thanks to the excision property of singular homology.

Moreover, $\varphi \in C^{1}(X)$ satisfies the $C$-condition if it has the following property:

"Every sequence $\left\{u_n\right\}_{n \in \mathbb{N}}\subset X$ such that $\left\{\varphi\left(u_n\right)\right\}_{n \in \mathbb{N}} \subset \mathbb{R}$   is bounded, and $\left(1+\left\|u_n\right\|_X\right) \varphi^{\prime}\left(u_n\right) \rightarrow 0$ in $X^*$, admits a strongly convergent subsequence."

Suppose $\varphi \in C^{1}(X)$ satisfies the $C$-condition and $\inf \varphi(K_{\varphi})>-\infty$.  We can define the critical groups of $\varphi(\cdot)$ at infinity when $c<\inf\varphi(K_{\varphi})$, denoted by
$$
C_k(\varphi, \infty)=H_k\left(X, \varphi^c\right) \,\,\text {for all} \,\,k \in \mathbb{N}_0.
$$
These critical groups are independent of the level $c<\inf \varphi(K_{\varphi})$ and are well-defined thanks to the second deformation theorem (see \cite[Theorem~5.3.12]{PRR}).

Suppose that $K_{\varphi}$ is finite, we introduce  the following series  in $t \in \mathbb{R}$,
$$
\begin{aligned}
& M(t, u)=\sum_{k \in \mathbb{N}_0} \operatorname{rank} C_k(\varphi, u) t^k \,\, \text {with}\,\, u\in K_{\varphi}, \\
& P(t, \infty)=\sum_{k \in \mathbb{N}_0} \operatorname{rank} C_k(\varphi, \infty) t^k .
\end{aligned}
$$
The "Morse relation" says that
\begin{equation}\label{001}
\sum_{u \in K_{\varphi}} M(t, u)=P(t, \infty)+(1+t) Q(t)
\end{equation}
with $Q(t)=\sum_{k \in \mathbb{N}_0} \beta_k t^k$ a formal series in $t \in \mathbb{R}$ with nonnegative integer coefficients.

To use the properties of critical groups, we require the following notion. Let us define that $\varphi: \Omega \times \mathbb{R} \rightarrow \mathbb{R}$ as an $L^{\infty}$-locally Lipschitz integrand if, for all $x \in \mathbb{R}$, $ z \rightarrow \varphi(z, x)$ is measurable and for every compact set $K \subseteq \mathbb{R}$, there exists $\vartheta_K \in L^{\infty}(\Omega)$ such that
 $$|\varphi(z, x)-\varphi(z, y)| \leq \vartheta_K(z)|x-y|, \,\, \text{for a.e.} \,\, z \in \Omega, \,\, \text{all} \,\,  x, y \in K.$$
Clearly, such a function is jointly measurable (see \cite[Proposition 2.2.31]{PW}). Therefore, if $u \in L^0(\Omega)$, then $z \rightarrow \varphi(z, u(z))$ is measurable.

Suppose $u \in L^{0}(\Omega)$ has the property that for all compact subset $\mathcal{K} \subseteq \Omega$, there exists a constant $C_{\mathcal{K}}$ such that
 $$0<C_\mathcal{K} \leqslant u(z) \,\, \text{for a.e.}\,\, z\in \mathcal{K}  ,$$
   then we denote that $0\prec u$. Similarly,  $v\prec 0$ is used if $0\prec -v$.

Let $\hat{\lambda}_1(q)$  denote the principal eigenvalue of $\left(-\Delta_q, W_0^{1, q}(\Omega)\right)$. We know that $\hat{\lambda}_1(q)>0$ and it is simple and isolated. It has the following variational characterization:
\begin{equation}\label{chara}
\hat{\lambda}_1(q)=\inf\left\{\frac{\|\nabla u\|_q^q}{\|u\|_q^q}: u \in W_0^{1, q}(\Omega), \,\,  u \neq 0\right\}.
\end{equation}
The infimum in \eqref{chara} is realized on the corresponding one-dimensional eigenspace,  the elements of which have a fixed sign. In fact, $\hat{\lambda}_1(q)>0$ is the only eigenvalue with eigenfunctions of constant sign, while all other eigenvalues have nodal eigenfunctions. Using these facts, we can easily prove the following result (see \cite[Lemma~4.11]{MP}).

\begin{proposition}\label{eigen}
If $\vartheta\in L^{\infty}(\Omega)$ satisfies $\vartheta(z) \leqslant \hat{\lambda}_1(q)$ for a.e. $z \in \Omega$ and the inequality is strict on a set of positive Lebesgue measure, then there exists $c_1>0$ such that
$$
c_1\|\nabla u\|_q^q \leqslant\|\nabla u\|_q^q-\int_{\Omega}\vartheta(z)\vert u\vert^q d z,\,\, \text{for all}\,\, u\in W_0^{1, q}(\Omega).
$$

\end{proposition}

The hypotheses on the perturbation $f(z, x)$ are the following:\\
$(H_1):  f: \Omega \times \mathbb{R} \rightarrow \mathbb{R}$ is an $L^{\infty}$-locally Lipschitz integrand such that $f(z, 0)=0$ for a.e. $z \in \Omega$ and

\begin{enumerate}
	\item $\vert f(z, x)\vert\leq \hat{\alpha}(z)(1+|x|^{r-1})$ for a.e. $z\in\Omega$, all $x\in \mathbb{R}$ with $\hat{\alpha} \in L^{\infty}(\Omega)$ and $p<r<q^*$;
\item there exist $\vartheta \in L^{\infty}(\Omega)$ and $\delta>0$ such that
$$
\begin{aligned}
& \vartheta(z) \leq \hat{\lambda}_1(q) \,\,\text { for a.e.} \,\, z \in \Omega, \,\,\vartheta \not\equiv \hat{\lambda}_1(q), \\
& \limsup _{x \rightarrow 0} \frac{q F(z, x)}{|x|^q} \leq\vartheta(z)\text { uniformly for a.e } z \in \Omega,
\end{aligned}
$$
where $F(z, x)=\int_0^x f(z, s) d s$ and $0 \leq f(z, x)x$   for a.e.  $z \in \Omega$, all $\vert x\vert\leq \delta$.
\end{enumerate}

\begin{remark}\label{Remarkper}
The hypotheses on the perturbation $f(z, x)$ are minimal, requiring a nonuniform nonresonance condition as $x \rightarrow 0$ and a local sign condition.
\end{remark}

Now let us consider the following auxiliary double phase Dirichlet problem
\begin{align}\label{111}
 \left\{
\begin{aligned}
&-\Delta_p^\alpha u(z)-\Delta_q u(z)=\lambda \vert u(z)\vert^{\tau-2}u(z), \,\,\,
\hbox{in }\Omega,\\
&u|_{\partial \Omega}=0,\,\, 1<\tau< q<p<N,\,\, \lambda>0.
\end{aligned}
\right.
\tag{$A_{\lambda}$}
\end{align}
By reasoning similarly to the proof of \cite[Proposition 10]{LP}, we have the following result concerning problem \eqref{111}.

\begin{proposition}\label{H0}
If hypotheses $(H_0)$ hold and $\lambda>0$, then problem \eqref{111} has a unique positive solution
$\bar{u}_\lambda \in W_0^{1, \eta}(\Omega) \cap L^{\infty}(\Omega)\,\,\text{with}\,\,0\prec \bar{u}_\lambda$. Furthermore, since problem \eqref{111} is odd,  $\bar{v}_\lambda=-\bar{u}_\lambda$ is the unique negative solution of \eqref{111}.
\end{proposition}

\section{Constant-sign solutions}\label{Sec3}
We define the following sets
$$S_\lambda^{+}= \text{set of positive solutions of problem $(P_{\lambda})$},$$
$$S_\lambda^{-}=\text{set of negative solutions of problem $(P_{\lambda})$}.$$

Now we will show that both $S_\lambda^{+}$ and $S_\lambda^{-}$ are non-empty for small $\lambda> 0$.
\begin{proposition}\label{nonempty}
Assuming hypotheses $(H_0)$, $(H_1)$ hold, for all $\lambda>0$ small, we have
$$
\begin{aligned}
& \emptyset\neq S_\lambda^{+} \subseteq W_0^{1, \eta}(\Omega) \cap L^\infty(\Omega), 0\prec u,\,\,  \text {for all} \,\,u\in S_\lambda^{+}, \\
& \emptyset\neq S_\lambda^{-} \subseteq W_0^{1, \eta}(\Omega) \cap L^\infty(\Omega), v\prec 0,\,\, \text {for all } \,\,v\in S_\lambda^{-}.
\end{aligned}
$$
\end{proposition}
\begin{proof}
Let $\varphi_\lambda^{+}: W_0^{1, \eta}(\Omega) \rightarrow \mathbb{R}$ be the $C^{1}$-functional defined by
$$
\varphi_\lambda^{+}(u)=\frac{1}{p}\rho_{\alpha}(\nabla u)+\frac{1}{q}\|\nabla u\|_q^q-\frac{\lambda}{\tau}\|u^{+}\|_\tau^\tau-\int_{\Omega}F(z, u^{+})d z,\,\,\text{ for all}\,\, u \in W_0^{1, \eta}(\Omega),
$$
where $\rho_{\alpha}(\nabla u)=\int_{\Omega} \alpha(z)\vert \nabla u\vert^p dz$. From hypotheses $(H_1)$(1), (2),  for any $\varepsilon>0$, there exists $c_1=c_1(\varepsilon)>0$ such that
$$
F(z, x) \leqslant \frac{1}{q}(\vartheta(z) +\varepsilon)\vert x\vert^q+c_1\vert x\vert^r,\,\,\text {for a.e.}\,\, z \in \Omega,\,\,\text {all}\,\,x \in \mathbb{R}.
$$
Then for all $u\in W_0^{1, \eta}(\Omega)$, we have
\begin{equation}\label{inequa11}
\varphi_\lambda^{+}(u)\geq \frac{1}{p}\rho_{a}(\nabla u) +\frac{1}{q}\left(\|\nabla u\|_q^q-\int_{\Omega}\vartheta(z)|u|^qdz-\varepsilon\|u\|_q^q\right)  -\frac{\lambda}{\tau}\|u\|_\tau^r-c_1\|u\|_r^r.
\end{equation}
By Proposition \ref{eigen}, we have from \eqref{chara} that
$$
\|\nabla u\|_q^q-\int_{\Omega} \vartheta(z)\vert u\vert^q d z-\varepsilon \Vert u\Vert_{q}^q  \geq \left(c_1-\frac{\varepsilon}{\hat{\lambda}_1(q)}\right)\Vert \nabla u\Vert_q^q.
$$
Choosing $\varepsilon\in \left(0, \hat{\lambda}_1(q) c_1\right)$, we obtain
\begin{equation}\label{inequa12}
\|\nabla u\|_q^q-\int_{\Omega} \vartheta(z)\|u\|^q d z-\varepsilon\|u\|_q^q \geqslant c_2\|\nabla u\|_q^q\,\,\text{for some}\,\,c_2>0.
\end{equation}
Returning to \eqref{inequa11} and using \eqref{inequa12} and the fact that
$$W_0^{1, \eta}(\Omega) \hookrightarrow L^\tau(\Omega), L^r(\Omega) \,\,\text{continuously}, $$
we have
$$
\begin{aligned}
\varphi_\lambda^{+}(u) \geq c_3 \rho_\eta(\nabla u)-c_4\left(\|u\|^\tau+\|u\|^r\right),\,\text {for some}\,\, c_3, c_4>0.
\end{aligned}
$$
If $\Vert u\Vert\leq 1$, from Proposition \ref{normre}, one has
\begin{equation}\label{inequa13}
\varphi_\lambda^{+}(u) \geq\Big(c_3-c_4\left(\|u\|^{\tau-p}+\|u\|^{r-p}\right)\Big)\|u\|^p.
\end{equation}
We consider the function
$$
\gamma_\lambda(t)=\lambda t^{\tau-p}+t^{r-p}, \,\, t>0.
$$
Since $1<\tau<q<p$, it is easy to see that
$$
\gamma_\lambda(t) \rightarrow+\infty \text { as } t \rightarrow 0^{+}\,\, \text {and as } \,\,t \rightarrow+\infty \text {. }
$$
So, there exits $t_0>0$ such that
$$
\begin{aligned}
\gamma_\lambda\left(t_0\right)=\min_{t>0} \gamma_\lambda(t)& \Rightarrow \gamma_\lambda^{\prime}\left(t_0\right)=0\\
& \Rightarrow \lambda(p-\tau) t_0^{\tau-p}=(r-p) t_0^{r-p} \\
& \Rightarrow t_0(\lambda)=\left(\frac{\lambda(p-\tau)}{r-p}\right)^{\frac{1}{r-\tau}}.
\end{aligned}
$$
Evidently $t_0(\lambda) \rightarrow 0^{+}$ as $\lambda \rightarrow 0^{+}$ and then since $p-\tau<r-\tau$,
we have
$$
\gamma_\lambda\left(t_0(\lambda)\right) \rightarrow 0^{+}\,\, \text {as } \,\,\lambda \rightarrow 0^{+}.
$$
Hence we can find $\lambda_*>0$ such that
$$t_0(\lambda)<1\,\,\, \text{and}\,\,\, \gamma_\lambda\left(t_0(\lambda)\right)<\frac{c_3}{c_4}\,\,\text{ for all} \,\,\lambda \in\left(0, \lambda_*\right).$$
Then from \eqref{inequa13}, we have that
\begin{equation}\label{inequa14}
\varphi_{\lambda}^{+}(u) \geq m_{\lambda}>0,\,\, \text{for all}\,\, \Vert u\Vert=t_0(\lambda)\,\,\text{and all}\,\, \lambda \in\left(0, \lambda_*\right).
\end{equation}

Now we introduce the following closed ball in $W_0^{1, \eta}(\Omega)$
$$
\bar{B}_\lambda=\left\{u \in W_0^{1, \eta}(\Omega): \Vert u \Vert\leq t_0(\lambda), \lambda \in\left(0, \lambda_*\right) \right\}.
$$
Since $W_0^{1, \eta}(\Omega)$ is reflexive, by  the James and Eberlein-Smulian theorems (see \cite{PW}), we have that $\bar{B}_\lambda$ is sequentially $\omega$-compact. Also using Proposition \ref{embeddings}, we see that $\varphi_{\lambda}^{+}(\cdot)$ is sequentially weakly lower semicontinuous. So, by the Weierstrass-Tonelli theorem, there exists a $u_\lambda\in W_0^{1, \eta}(\Omega)$ such that
\begin{equation}\label{inequa15}
\varphi_\lambda^{+}\left(u_\lambda\right)=\inf \left\{\varphi_\lambda^+(u): u \in \bar{B}_\lambda\right\}.
\end{equation}

Let $u \in C_0^{1}(\bar{\Omega})$ with $u(z)>0$ for all $z \in \Omega$. Choose $t\in (0, 1)$ small enough such that
$$
tu\in  \bar{B}_\lambda \,\,\text { and} \,\, 0\leq tu(z) \leq\delta \,\,\text {for all}\,\, z \in\bar{\Omega},
$$
with $\delta>0$ as postulated by hypothesis $(H_1)(2)$. Then from the local sign condition $\left(\sec H_1(2)\right)$, recall $t\in (0, 1)$ and $q<p$, we have
$$
\begin{aligned}
\varphi_\lambda^{+}(t u) \leqslant \frac{t^q}{q} \rho_\eta(\nabla u)  -\frac{\lambda t^\tau}{\tau}\|u\|_\tau^\tau.
\end{aligned}
$$
Since $\tau<q$, choose $t\in (0,1)$ even smaller if necessary, we have from \eqref{inequa15} that
\begin{align}\label{inequa16}
& \varphi_\lambda^{+}(t u)<0,\nonumber \\
\Rightarrow & \varphi_\lambda^{+}\left(v_\lambda\right)<0=\varphi_\lambda^{+}(0),\nonumber \\
\Rightarrow & v_\lambda \neq 0.
\end{align}
From \eqref{inequa14} and \eqref{inequa16} it follows that
$$
0<\left\|v_\lambda\right\|<t_0(\lambda).
$$
Then from \eqref{inequa15}, we have  for all $h \in W_0^{1,\eta}(\Omega)$ and $\lambda \in\left(0, \lambda_*\right)$ that
\begin{equation}\label{inequa17}
\left\langle\left(\varphi_\lambda^{+}\right)^{\prime}\left(u_\lambda\right), h\right\rangle=0\Rightarrow\left\langle V\left(u_\lambda\right), h\right\rangle=\int_{\Omega} \lambda\left(u_\lambda^{+}\right)^{\tau-1} h d z+\int_{\Omega} f\left(z, u_\lambda^{+}\right) h d z.
\end{equation}
In \eqref{inequa17}, choose the test function $h=-u^{-}_\lambda \in  W_0^{1, \eta}(\Omega)$ and use Proposition \ref{normre}, we obtain
$$
\rho_\eta\left(\nabla u_\lambda^-\right)=0 \Rightarrow u_\lambda \geqslant 0,\,\, u_\lambda \neq 0.
$$
Then from \eqref{inequa17} we see that
$$
u_\lambda\in S_\lambda^{+},\,\,\text {for all}\,\, \lambda \in\left(0, \lambda_*\right).
$$
\cite[Theorem 3.1]{GW} implies that
$$
u_\lambda \in W_0^{1,\eta}(\Omega) \cap L^{\infty}(\Omega).
$$
Let $\rho=\left\|u_\lambda\right\|_{\infty}$. Hypotheses $(H_1)$ imply that we can find $\hat{\xi}_\rho>0$ such that
$$
f\left(z, x\right)+\hat{\xi}_\rho x^{p-1}\geq 0,\,\,\text {for a.e.} \,\,z \in \Omega, \text { and all}\,\, 0\leq x\leq \rho.
$$
Therefore we have
$$
-\Delta_p^\alpha u_{\lambda}-\Delta_q u_\lambda+\hat{\varepsilon}_p u_\lambda^{p-1} \geq 0\,\, \text { in }\, \Omega.
$$
Then using \cite[Proposition 2.4]{PVV} we obtain
$$
0\prec u_\lambda.
$$
Therefore we have proved that for $\lambda\in \left(0, \lambda_*\right)$
$$
\emptyset\neq S_\lambda\subseteq W_0^{1,\eta}(\Omega) \cap L^{\infty}(\Omega) \text { and } 0\prec u \,\,\text {for all}\,\, u \in S_\lambda^{+}.
$$

Similarly we can show the nonemptiness of the set $S_\lambda^{-}$. In this case we work with the $C^{1}$-functional $\varphi_\lambda^{-}: W_0^{1, \eta}(\Omega) \rightarrow \mathbb{R}$  defined by
$$
\varphi_\lambda^{-}(u)=\frac{1}{p}\rho_{a}(\nabla u)+\frac{1}{q}\|\nabla u\|_q^q-\frac{\lambda}{\tau}\|u^{-}\|_\tau^\tau-\int_{\Omega}F(z, u^{-})d z,\,\,\text{ for all}\,\, u \in W_0^{1, \eta}(\Omega).
$$
\end{proof}
The proof of the previous proposition, leads to the next result which will be used to show that the nodal solutions that we will produce, asymptotically vanish as $\lambda \rightarrow 0^{+}$.
\begin{proposition}\label{vanishing}
 If hypotheses $(H_0), (H_1)$ hold and $\lambda\in \left(0, \lambda_*\right)$, then there exist $u_{\lambda} \in S_\lambda^{+}$ and $v_\lambda \in S_\lambda^{-}$ such that
$$
u_{\lambda}, v_{\lambda}\rightarrow 0,\,\, \text {in}\,\, W_0^{1,\eta}(\Omega) \cap L^{\infty}(\Omega)\,\, \text {as} \,\,\lambda \rightarrow 0^{+}.
$$
\end{proposition}

\begin{proof} From the proof of Proposition \ref{nonempty}, we know that we can find $u_{\lambda} \in S_\lambda^{+}$ such that
$$
\Vert u_\lambda\Vert<t_0(\lambda).
$$
Since $t_0(\cdot)$ is increasing and $t_0(\lambda) \rightarrow 0$ as $\lambda \rightarrow 0^{+}$, from \cite[Theorem 3.1]{GW} (see also \cite[Theorem 3.1]{MW}, we infer that
$$
\{u_{\lambda}\}_{\lambda\in \left(0, \lambda_*\right)}\subseteq  L^{\infty}(\Omega)\,\,\text{is bounded and}\,\, \Vert u_{\lambda}\Vert_{\infty}\leq O(\lambda),\,\,\lambda\in \left(0, \lambda_*\right),
$$
here $O(\lambda)\rightarrow 0$ as $\lambda \rightarrow 0^{+}$. Therefore finally we have
$$
u_\lambda \rightarrow 0 \,\,\text { in } \,\, W_0^{1,\eta}(\Omega) \cap L^{\infty}(\Omega)\,\, \text { as }\,\, \lambda \rightarrow 0^+.
$$

Similarly, we have solutions $v_\lambda \in S_\lambda^{-}$ for all $ \lambda \in\left(0, \lambda_*\right)$ such that
$$
v_\lambda \rightarrow 0 \,\,\text { in } \,\, W_0^{1,\eta}(\Omega) \cap L^{\infty}(\Omega)\,\, \text { as }\,\, \lambda \rightarrow 0^+.
$$
\end{proof}

Next we will show that for $\lambda \in\left(0, \lambda_*\right)$, problem $(P_{\lambda})$ has extremal constant-sign solutions (that is, a largest positive solution and a smallest negative solution). We will use them to generate a nodal (sign-changing) solution.

So, we have the following result.

\begin{proposition}\label{sign}
If hypotheses $(\mathrm{H}_0), (\mathrm{H}_{1})$ hold and $ \lambda \in\left(0, \lambda_*\right)$,
then there exist $u_\lambda^* \in S_\lambda^{+}$ and $v_\lambda^*\in S_\lambda^{-}$ such that
$$
\begin{aligned}
u_\lambda^* \leqslant u,\,\,\text{ for all}\,\, u \in S_\lambda^{+},\\
v \leqslant v_\lambda^*,\,\,\text{ for all}\,\, v \in S_\lambda^{-}.
\end{aligned}
$$
\end{proposition}
\begin{proof}
As demonstrated in \cite{FP}, we show that
$$
S_\lambda^{+} \text { is downward directed. }
$$
Consequently, by using \cite[Theorem 5.109]{HP}, we can obtain a decreasing sequence $\{u_n \}_{n\in \mathbb{N}} \subseteq S_\lambda^{+}$ such that
 $$
\inf S_\lambda^{+}=\inf_{n \in \mathbb{N}} u_n.
$$
We have
\begin{equation}\label{inequa21}
\left\langle V(u_n),  h\right\rangle=\int_{\Omega} \lambda u_n^{\tau-1} h d z+\int_{\Omega}f\left(z, u_n\right) h dz, \,\, \text{for all}\,\,h \in W_0^{1,\eta}(\Omega),\,\, \text{all}\,\, n \in \mathbb{N},
\end{equation}
and
\begin{equation}\label{inequa22}
0 \leq u_n \leq u_1,\,\,\text{ for all}\,\, n \in \mathbb{N}.
\end{equation}
Based on hypotheses  $(H_1)$, we have refer from \eqref{inequa22} and the fact that $ u_1 \in L^{\infty}(\Omega)$ that
$$
\begin{aligned}
\left|\lambda u_n(z)^{\tau-1}+f\left(z, u_n(z)\right)\right| &\leq c_5\left(u_n(z)^{\tau-1}+u_n(z)^{q-1}+u_n(z)^{r-1}\right)\\
& \leq c_5\left(u_1(z)^{\tau-1}+u_1(z)^{q-1}+u_1(z)^{r-1}\right) \\
& \leq c_6,\,\,\text{for a.e.}\,\, z\in \Omega,
\end{aligned}
$$
where $c_5$ and $c_6$ are positive constants. Then using Moser's iteration technique, as in \cite[Proposition 1.3]{GV}, we have
\begin{equation}\label{inequa23}
\left\|u_n\right\| \leqslant c_7 O\left(\left\|f\left(\cdot, u_n(\cdot)\right)\right\|_m\right), \,\,\text{for some}\,\, c_7>0, \,\,\text{with}\,\, m>N,\,\, \text{all}\,\, n \in \mathbb{N}.
\end{equation}
 From \eqref{inequa21}, using the test function $h=u_n\in W_0^{1, \eta}(\Omega)$, we obtain that
$$
\left\{u_n\right\}_{n \in \mathbb{N}} \subseteq W_0^{1, \eta}(\Omega)\,\, \text { is bounded }.
$$
So, from Proposition \ref{embeddings} we may assume that
\begin{equation}\label{inequa24}
u_n \rightharpoonup u_\lambda^* \,\,\text {in}\,\, W_0^{1, \eta}(\Omega), \,\, u_n \rightarrow u_\lambda^* \,\,\text {in} \,\,L^p(\Omega).
\end{equation}
In \eqref{inequa21} we use the test function $h=(u_n - u_\lambda^*) \in W_0^{1, \eta}(\Omega)$, pass to the limit as $n \rightarrow \infty$ and use \eqref{inequa24}. Then, by Proposition \ref{mappro}, we have
\begin{equation}\label{inequa31}
\lim _{n \rightarrow \infty}\left\langle V\left(u_n\right), u_n-u_\lambda^*\right\rangle=0\Rightarrow  u_n \rightarrow u_\lambda^*\,\,\text{in}\,\,W_0^{1, \eta}(\Omega).
\end{equation}
Suppose that $u_\lambda^*=0$. Then from \eqref{inequa23} it follows that
$$
u_n \rightarrow 0 \,\,\text {in} \,\, L^{\infty}(\Omega)\,\, \text { as }\,\, n \rightarrow +\infty.
$$
Therefore there exists $n_0 \in \mathbb{N}$ from hypothesis $(H_1)(2)$ such that for a.e. $z \in \Omega$ and all $n \geqslant n_0$
\begin{equation}\label{inequa25}
0 \leq u_n(z) \leq\delta \Rightarrow \lambda u_n(z)^{\tau-1} \leqslant \lambda u_n(z)^{\tau-1}+f\left(z, u_n(z)\right).
\end{equation}

We fix $n \geqslant n_0$ and introduce Carath\'{e}odory function
$\gamma_\lambda^{+}(z, x)$ defined by
\begin{equation}\label{inequa26}
\gamma_\lambda^{+}(z, x)= \begin{cases}\lambda\left(x^{+}\right)^{\tau-1} & \text { if } x \leqslant u_n(z), \\ \lambda u_n(z)^{\tau-1} & \text { if } u_n(z)<x .\end{cases}
\end{equation}
Assume $\Gamma_\lambda^{+}(z, x)=\int_0^t \gamma_\lambda^{+}(z, s) d s$ and consider the $C^{1}$-functional $\sigma_\lambda^{+}: W_0^{1, \eta}(\Omega) \rightarrow \mathbb{R}$ defined by
$$
\sigma_\lambda^{+}(u)=\frac{1}{p}\rho_{a}(\nabla u)+\frac{1}{q}\|\nabla u\|_q^q-\int_{\Omega}\Gamma_\lambda^{+}(z, u)d z,\,\,\text{ for all}\,\, u \in W_0^{1, \eta}(\Omega).
$$

From \eqref{inequa26} we see that $\sigma_\lambda^{+}(\cdot)$ is coercive. Also, it is sequentially weakly lower semicontinuous. So, we can find $\tilde{u}_\lambda \in W_0^{1, \eta}(\Omega)$ such that
\begin{equation}\label{inequa27}
\sigma_\lambda^{+}\left(\tilde{u}_\lambda\right)=\inf\left\{\sigma_\lambda^{+}(u): u \in W_0^{1, \eta}(\Omega)\right\} .
\end{equation}

Let $u \in C_0^1(\bar{\Omega})$ with $u(z)>0$ for all $z \in \Omega$ and let $t\in (0,1)$, then by \eqref{inequa26}, we have
$$
\begin{aligned}
\sigma_\lambda^{+}(t u) \leq &\frac{t^q}{q} \rho_\eta(\nabla u)-\int_{\Omega} \Gamma_\lambda^{+}(z, t u) d z \\
 =&\frac{t^q}{q} \rho_\eta(\nabla u)-\frac{\lambda t^\tau}{\tau} \int_{\left\{0 \leq t u \leq u_n\right\}} u^\tau d z \\
& -\frac{\lambda}{\tau} \int_{\left\{u_n<t u\}\right.} u_n^\tau d z-\lambda \int_{\left\{u_n<t u\}\right.} u_n^\tau\left(t u-u_n\right) d z \\
\leq & \frac{t^q}{q} \rho_\eta(\nabla u)-\frac{\lambda t^\tau}{\tau} \int_{\left\{0 \leq t u \leq u_n\right\}} u^\tau d z \\
 =&\frac{t^q}{q} \rho_\eta(\nabla u)-\frac{\lambda t^\tau}{\tau} \int_{\Omega} u^\tau d z+\frac{\lambda t^\tau}{\tau} \int_{\left\{u_n< t u\}\right.} u^\tau dz.
\end{aligned}
$$
Moreover, we have
\begin{equation}\label{inequa28}
\frac{\sigma_\lambda^{+}(t u)}{t^\tau} \leqslant \frac{t^{q-\tau}}{q} \rho_\eta(\nabla u)-\frac{\lambda}{\tau}\|u\|_\tau^\tau+\frac{\lambda}{\tau} \int_{\left\{u_n<t u\}\right.} u^\tau d z.
\end{equation}
Note that
$$
\begin{aligned}
& \frac{t^{q-\tau}}{q} \rho_\eta(\nabla u) \rightarrow 0\,\, \text { as }\,\, t \rightarrow 0^{+}, \\
& \frac{\lambda}{\tau} \int_{\left\{u_n<t u\right\}} u^\tau d z \rightarrow 0\,\, \text { as }\,\, t \rightarrow 0^{+}.
\end{aligned}
$$
Therefore from \eqref{inequa28} we have
$$
\limsup _{t \rightarrow 0^{+}} \frac{\sigma_\lambda^{+}(t u)}{t^\tau}=-\frac{\lambda}{\tau}\|u\|_\tau^\tau<0 .
$$
So, for $t \in(0,1)$ small, one has
\begin{align*}
& \sigma_\lambda^{+}(t u)<0, \\
\Rightarrow & \sigma_\lambda^{+}\left(\tilde{u}_\lambda\right)<0=\sigma_\lambda^{+}(0) \\
\Rightarrow & \tilde{u}_\lambda \neq 0 .
\end{align*}
From \eqref{inequa27} we have  for all $h\in W_0^{1, \eta}(\Omega)$ that
$$
\begin{aligned}
\left\langle\left(\sigma_\lambda^{+}\right)^{\prime}\left(\tilde{u}_\lambda\right), h\right\rangle=0 \,\,\Rightarrow  \left\langle V\left(\tilde{u}_\lambda\right), h\right\rangle=\int_{\Omega} \gamma_\lambda^{+}\left(z, \tilde{u}_\lambda\right) hdz.
\end{aligned}
$$
Choosing $h=-\tilde{u}_\lambda^{-} \in W_0^{1, \eta}(\Omega)$, by \eqref{inequa26}, we obtain
$$
\begin{aligned}
 \rho_\eta(\nabla \tilde{u}_\lambda^{-})=0 \Rightarrow  \tilde{u}_\lambda\geq 0,\,\, \tilde{u}_\lambda \neq 0.
\end{aligned}
$$
Also, if we use the test function $h=\left(\tilde{u}_\lambda-u_n\right)^{+}\in W_0^{1, \eta}(\Omega)$, since $u_n \in S_\lambda^{+}$ and \eqref{inequa25}, we have
$$
\begin{aligned}
& \left\langle V\left(\tilde{u}_\lambda\right),\left(\tilde{u}_\lambda-u_n\right)^{+}\right\rangle \\
=&\int_{\Omega} \lambda u_n^{\tau-1}\left(\tilde{u}_\lambda-u_n\right)^{+}d z \\
\leq & \int_{\Omega}\left(\lambda u_n^{\tau-1}+f\left(z, u_n\right)\right)\left(\tilde{u}_\lambda-u_n\right)^{+} dz \\
=& \left\langle V\left(u_n\right),\left(\tilde{u}_\lambda-u_n\right)^{+}\right\rangle \text , \\
\Rightarrow & \tilde{u}_\lambda\leq u_n.
\end{aligned}
$$
So, we have proved that
$$
\tilde{u}_\lambda \in\left[0, u_n\right],\,\, \tilde{u}_\lambda \neq 0.
$$
From \eqref{inequa26} it follows that $\tilde{u}_\lambda$ is a positive solution of \eqref{111}.
Then by Proposition \ref{H0}, we have
$$
\begin{aligned}
\tilde{u}_\lambda=\bar{u}_\lambda \Rightarrow  \bar{u}_\lambda\leq u_n \,\,\text { for all }\,\, n \geqslant n_0.
\end{aligned}
$$
But this contradicts the hypothesis that $u_n \rightarrow 0$ in
$W_0^{1, \eta}(\Omega)$ (recall that we have assumed that $u_\lambda^*=0$). Therefore $u_\lambda^* \neq 0$. In \eqref{inequa21}, pass to the limit as $n \rightarrow \infty$ and use \eqref{inequa31}, we obtain
$$
\left\langle V\left(u_\lambda^*\right), h\right\rangle=\int_{\Omega} \lambda\left(u_\lambda^*\right)^{\tau-1} h d z+\int_{\Omega} f\left(z, u_\lambda^*\right) h d z,\,\, \text { for all }\,\, h \in W_0^{1, \eta}(\Omega).
$$
Moreover, we have
$$
u_\lambda^* \in S_\lambda^{+},\,\, u_\lambda^*=\inf S_\lambda^{+}.
$$

Similarly we produce a maximal element for $S_\lambda^{-}$. The set $S_\lambda^{-}$ is upward directed. So, we can find $\left\{v_n\right\}_{n \in \mathbb{N}}\subseteq S_\lambda^{-}$ increasing sequence such that $\sup S_\lambda^{-}=\sup _{n \in \mathbb{N}} v_n$.
\end{proof}

In the next section, we will use $u_\lambda^*$ and $v_\lambda^*$ to generate a nodal solution. The idea is to  look for nontrivial solutions of $(P_{\lambda})$ $(\lambda\in (0, \lambda_*))$ in the order interval $\left[v_\lambda^*, \, u_\lambda^*\right]$ which are different from $u_\lambda^*$ and $v_\lambda^*$. On account of the extremality of $u_\lambda^*$ and $v_\lambda^*$, such a solution must be nodal. To produce this solution, we shall use truncation and comparison techniques and critical groups.

\section{Nodal solutions}\label{Sec4}

In order to focus on the order interval $\left[v_\lambda^*, u_\lambda^*\right]$ $\left(\lambda\in \left(0,  \lambda_*\right)\right)$, we introduce the following truncation of the reaction of the reaction of $(P_{\lambda})$
\begin{align}\label{inequa41}
 g_\lambda(z, x)= \begin{cases}\lambda\vert v_\lambda^*(z)\vert^{\tau-2} v_\lambda^*(z)+f\left(z, v_\lambda^*(z)\right) & \text { if } x<v_\lambda^*(z), \\ \lambda\vert x\vert^{\tau-2} x+f(z, x) & \text { if } v_\lambda^*(z) \leq x\leq  u_\lambda^*(z), \\ \lambda u_\lambda^*(z)^{\tau-1}+f\left(z, u_\lambda^*(z)\right) & \text { if } u_\lambda^*(z)<x.\end{cases}
\end{align}
Also we consider the positive and negative truncations
of $g_\lambda(z, \cdot)$, namely the functions
\begin{equation}\label{inequa42}
g_\lambda^{ \pm}(z, x)=g_\lambda(z, \pm x^{\pm}).
\end{equation}
All three are Carath\'{e}odory functions, we set
$$
G_\lambda(z, x)=\int_0^x g_\lambda(z, s) d s \text { and } G_\lambda^\pm(z, x)=\int_0^x g_\lambda^{ \pm}(z, s) d s,
$$
and consider the $C^{1}$-functionals $\beta_\lambda$, $\beta_\lambda^{ \pm}: W_0^{1, \eta}(\Omega) \rightarrow \mathbb{R}$ defined by
$$
\begin{aligned}
& \beta_\lambda(u)=\frac{1}{p}\rho_{a}(\nabla u)+\frac{1}{q}\|\nabla u\|_q^q-\int_{\Omega} G_\lambda(z, u) d z, \\
& \beta_\lambda^{ \pm}(u)=\frac{1}{p}\rho_{a}(\nabla u)+\frac{1}{q}\|\nabla u\|_q^q-\int_{\Omega} G_\lambda^\pm(z, u) d z,
\end{aligned}
$$
for all $u \in W_0^{1, \eta}(\Omega)$. From \eqref{inequa41} and \eqref{inequa42}, we can see that
$$
K_{\beta_\lambda} \subseteq \left[v_\lambda^*, u_\lambda^*\right],\,\, K_{\beta_\lambda^+} \subseteq\left[0, u_\lambda^*\right],\,\, K_{\beta_\lambda^-} \subseteq\left[v_{\lambda}^*,  0\right].
$$
The extremality of $u_\lambda^*, v_\lambda^*$ implies that
\begin{equation}\label{inequa43}
K_{\beta_\lambda} \subseteq \left[v_\lambda^*,\,\, \,u_\lambda^*\right], K_{\beta_\lambda^+} =\{0, u_\lambda^*\}, \,\,K_{\beta_\lambda^-}=\{v_{\lambda}^*,  0\}.
\end{equation}

Also let $\varphi_\lambda: W_0^{1, \eta}(\Omega) \rightarrow \mathbb{R}$ be the energy functional of problem $(P_{\lambda})$ defined by
$$
\varphi_\lambda(u)=\frac{1}{p} \rho_\alpha(\nabla u)+\frac{1}{q}\|\nabla u\|_q^q-\frac{\lambda}{\tau}\|u\|_\tau^\tau-\int_{\Omega} F\left(z, u\right) d z, \,\, \text { for all }\,\, u\in W_0^{1, \eta}(\Omega).
$$
Evidently $\varphi_\lambda \in C^{1}\left(W_0^{1, \eta}(\Omega)\right)$.

As mentioned in the introduction, we will address the challenges arising from the absence of a global regularity theory by employing critical groups. We will compute the critical groups of $\beta_\lambda(\cdot)$ and $\beta_\lambda^{ \pm}(\cdot)$.

First of all,  we will compute the critical groups of $\beta_\lambda(\cdot)$ at $0$.

\begin{proposition}\label{crigro}
If hypotheses $(H_{0}), (H_1)$ hold and $\lambda\in \left(0, \lambda_*\right)$,
then
$$
C_k\left(\beta_\lambda, 0\right)=0,\,\, \text{ for all }\,\, k \in \mathbb{N}_0.
$$
\end{proposition}
\begin{proof} For any $u\in W_0^{1, \eta}(\Omega)$, we have
\begin{align}\label{inequa44}
& \left|\varphi_\lambda(u)-\beta_\lambda(u)\right| \nonumber\\
=&\left|\int_{\Omega}\left(\frac{\lambda}{\tau}|u|^\tau+F(z, u) - G_\lambda(z, u)\right) dz\right| \nonumber\\
\leq &\int_{\left\{u<v_\lambda^*\right\}}\left|\frac{\lambda}{\tau}(|u|^\tau-|v_\lambda^*|^r )-\lambda\vert v_\lambda^*\vert^{\tau-2} v_\lambda^*\left(u-v_\lambda^*\right) \right| d z \nonumber\\
& + \int_{\left\{u<v_\lambda^*\right\}}\left\vert F(z, u)-F\left(z, v_\lambda^*\right)-f\left(z, v_\lambda^*\right)\left(u-v_\lambda^*\right) \right\vert d z \nonumber\\
& +\int_{\left\{u_\lambda^*<u\right\}}\left|\frac{\lambda}{\tau}\Big(u^\tau-(u_\lambda^*)^r \Big)-\lambda(u_\lambda^*)^{\tau-1} \left(u-u_\lambda^*\right) \right| d z \nonumber\\
& + \int_{\left\{u_\lambda^*<u\right\}}\left\vert F(z, u)-F\left(z, u_\lambda^*\right)-f\left(z, u_\lambda^*\right)\left(u-u_\lambda^*\right) \right\vert dz.
\end{align}
Note that, from the continuous embedding $W_0^{l, \eta}(\Omega) \hookrightarrow L^\tau(\Omega)$,  we have that for some $c_8>0$
\begin{equation}\label{inequa45}
\int_{\left\{u<v_\lambda^*\right\}}\left|\frac{\lambda}{\tau}(|u|^\tau-|v_\lambda^*|^r )-\lambda\vert v_\lambda^*\vert^{\tau-2} v_\lambda^*\left(u-v_\lambda^*\right) \right| d z\leq \lambda c_8\|u\|^\tau.
\end{equation}
Since $F(z, \cdot)$ is $L^{\infty}$-locally Lipschitz integrand and
$$
\left|v_\lambda^*\right| \leqslant|u|,\,\, \text { on }\,\,\left\{u<v_\lambda^*\right\},
$$
it yields
$$\int_{\left\{u<v_\lambda^*\right\}}\left|F(z, u)-F\left(z, v_\lambda^*\right)-f\left(z, v_\lambda^*\right)\left(u-v_\lambda^*\right) \right| d z\leq  c_9\|u\|
$$
for some $c_9>0$.

Similarly, we have for some $c_{10}>0$
\begin{equation}\label{inequa46}
	\int_{\{u^*_\lambda<u\}}\left|\frac{\lambda}{\tau}(u^\tau-(u^*_\lambda)^\tau)-\lambda(u^*_\lambda)^{\tau-1}(u-u^*_\lambda)\right|dz\leq\lambda c_{10}\|u\|^\tau
\end{equation}
and for some $c_{11}>0$
\begin{equation}\label{inequa47}
	\int_{\{u^*_\lambda<u\}}\left|F(z,u)-F(z,u^*_\lambda)-f(z,u^*_\lambda)(u-u^*_\lambda)\right|dz\leq c_{11}\|u\|.
\end{equation}

Returning to \eqref{inequa44} and using \eqref{inequa45}, \eqref{inequa46} and \eqref{inequa47}, we have for some $c_{12} >0$ and all $\|u\|\leq1$
\begin{equation}\label{inequa48}
	\left| \varphi _\lambda(u)-\beta_\lambda(u) \right|\leq c_{12}\Vert u\Vert^\tau.
\end{equation}

Next we conduct a similar estimation for the difference of the two derivatives. So let $u,h\in W^{1,\eta}_0(\Omega)$, we have
\begin{align*}
	& \left|\left\langle\varphi_{\lambda}'(u)-\beta_{\lambda}^{\prime}(u), h\right\rangle\right| \\
	\leq & \int_{\left\{u<v_{\lambda}^{*}\right\}}\lambda\left||u|^{\tau-2} u-\left|v_{\lambda}^{*}\right|^{\tau-2} v_{\lambda}^{*}\right| |h|dz+\int_{\left\{u< v_{\lambda}^{*}\right\}}\left| f(z, u) - f\left(z, v_{\lambda}^{*}\right)\right||h| d z \\
	&+ \int_{\{u_{\lambda}^{*}<u\}} \lambda\left|u^{\tau-1}-\left(u_{\lambda}^{*}\right)^{\tau-1}\right||h| dz+\int_{\{u_{\lambda}^{*}<u\}} \left| f(z, u) -f(z,u_{\lambda}^{*})\right|\vert h\vert d z.
\end{align*}

Note that $|u|^{\tau-1}\in L^{\tau'}(\Omega)~\left(\frac{1}{\tau}+\frac{1}{\tau'}=1\right)$, while from Proposition \ref{embeddings}, we know that $\vert h\vert\in L^\tau(\Omega)$.  Also $\vert u\vert\in L^{(q^*)'}(\Omega), \vert h\vert\in L^{q^*}(\Omega)$, here $2\leq q^*$. So, using the H\"older's inequality and Proposition \ref{embeddings}, we have for some $c_{13}>0$ and all $\vert u\vert\leq 1$,
\begin{equation}\label{inequa49}
\left|\left\langle\varphi_{\lambda}^{\prime}(u)-\beta_{\lambda}^{\prime}(u), h\right\rangle\right| \leq(\lambda+1) c_{13}\|u\|\|h\| \Rightarrow \|\varphi_{\lambda}^{\prime}(u)-\beta_{\lambda}^{\prime}(u)\|_*\leq(\lambda+1)c_{13}\Vert u\Vert.
\end{equation}

From \eqref{inequa48} and \eqref{inequa49}, for any $\varepsilon>0$, there exists $\delta\in(0,1)$ such that
\begin{align*}
	\|\varphi_\lambda-\beta_\lambda\|_{C^1(\bar B_\delta)}\leq\varepsilon
\end{align*}
where $\bar B_\delta=\{u\in W^{1,\eta}_0(\Omega):\|u\|\leq\delta\}$.

We assume that $K_{\beta_\lambda}$ is finite. Otherwise on account of \eqref{inequa43} we already have an infinity of nodal solutions and so we are done. Therefore we can use the $C^1$-continuity property of critical groups (see \cite[Theorem 5.126]{GP}) and obtain
\begin{equation}\label{inequa411}
	C_k(\varphi_\lambda,0)=C_k(\beta_\lambda,0) \quad\text{for~all}~ k\in\mathbb{N}_0.
\end{equation}

Since $0\leq\lambda\vert x\vert^\tau\leq\vert x\vert^\tau+f(z,x)x$ for a.e. $z \in \Omega$ and all $|x|\leq\delta$, then from  \cite[Proposition 9]{LP} and \eqref{inequa411}, we have
\begin{align*}
	C_k(\varphi_\lambda, 0)=0 \,\,\text{for~all}\,\, k\in\mathbb{N}_0 \Rightarrow C_k(\beta_\lambda, 0)=0 \,\,\text{for~all}\,\,k\in\mathbb{N}_0.
\end{align*}

\end{proof}
\begin{proposition}\label{pro10}
	If hypotheses $(H_0), (H_1)$ hold and $\lambda\in(0,\lambda_*)$, then
$$C_k(\beta_\lambda,u^*_\lambda)=C_k(\beta_\lambda^+,v^*_\lambda)\,\,\,\text{and}\,\,\,C_k(\beta_\lambda,v^*_\lambda)=C_k(\beta^-_\lambda,v^*_\lambda)\,\,\,\text{for all}\,\,\,k\in\mathbb{N}_0.
$$
\end{proposition}
\begin{proof}
Let $W_+=\left\{u\in W^{1,\eta}_0(\Omega): 0\leq u(z)\,\,\text{for a.e.}\,\, z\in\Omega\right\}$. Note that
\begin{equation}\label{inequa412}
	\beta_\lambda|_{W_+}=\beta_\lambda^{+}|_{W_+}.
\end{equation}
For any $u\in W^{1,\eta}_0(\Omega)$, by  \eqref{inequa41} and \eqref{inequa412} we have
\begin{align}\label{inequa413}
&\vert\beta_\lambda(u)-\beta^+_\lambda(u)\vert\notag\\
\leq &\int_{\Omega}|G_\lambda(z,u)-G^+_\lambda(z,u)|dz\notag\\
\leq &\int_\Omega|G_\lambda(z,u)-G_\lambda(z,u^*_\lambda)|dz+\int_\Omega|G^+_\lambda(z,u^*_\lambda)-G^+_\lambda(z,u)|dz.
\end{align}

We will now estimate the two integrals on the right-hand side of \eqref{inequa413},
\begin{align}\label{inequa414}
	&\int_{\Omega}|G_\lambda(z,u)-G_\lambda(z,u^*_\lambda)|dz\notag\\
	\leq &\int_{\{u<v^*_\lambda\}}\left| \frac{\lambda}{\tau}\left(|v^*_\lambda|^\tau-(u^*_\lambda)^\tau\right) +|v^*_\lambda|^{\tau-2}v^*_\lambda(u-v^*_\lambda)\right|dz\notag\\
	+&\int_{\{u<v^*_\lambda\}}\left|F(z,v^*_\lambda)-F(z,u^*_\lambda)+f(z,v^*_\lambda)(u-v^*_\lambda)\right|dz\notag\\
	+&\int_{\{v^*_\lambda\leq u\leq u^*_\lambda\}}\left| \frac{\lambda}{\tau}(|u|^\tau-(u^*_\lambda)^\tau)+F(z,u)-F(z,u^*_\lambda)\right|dz\notag\\
	+&\int_{\{u^*_\lambda<u\}}\left| \lambda(u^*_\lambda)^{\tau-1}(u-u^*_\lambda)+ f(z,u^*_\lambda)(u-u^*_\lambda)\right|dz.
\end{align}

We make the following observations.

$\bullet$  If $|v^*_\lambda|\leq u^*_\lambda$, then $\vert v^*_\lambda\vert^{\tau}-(u^*_\lambda)^\tau\leq 0$.

If $u^*_\lambda<|v^*_\lambda|$,  since $u^*_\lambda, v^*_\lambda\in L^\infty(\Omega)$, we have for some $c_{14}>0$
\begin{equation*}
	0 \leq\left|v_{\lambda}^{*}\right|^{\tau}-\left(u_{\lambda}^{*}\right)^{\tau} \leq c_{14}\left\{\begin{array}{l}
		\left(\left|v_{\lambda}^{*}\right|-u_{\lambda}^{*}\right)^{\tau} \quad\text{if}~ \tau\leq 2,\\
		\left|v_{\lambda}^{*}\right|-u_{\lambda}^{*}~\quad\quad\text{if}~ 2<\tau.
		\end{array}\right.
\end{equation*}

On $\{u<v^*_\lambda\}$, we have
\begin{align*}
	|v^*_\lambda|\leq|u|.
\end{align*}
Therefore on $\{u<v^*_\lambda\}$, we have
\begin{equation*}
\begin{aligned}
	0 \leq\left|v_{\lambda}^{*}\right|^{\tau}-\left(u_{\lambda}^{*}\right)^{\tau} & \leq c_{14}\left\{\begin{array}{l}
	\left(|u|-u_{\lambda}^{*}\right)^{\tau}\quad\text{if}~ \tau\leq 2,\\
	|u|-u_{\lambda}^{*}~\quad\quad\text{if}~ 2<\tau.
	\end{array}\right. \\
	& \leq c_{14}\left\{\begin{array}{l}
	\left|u-u_{\lambda}^{*}\right|^{\tau}\quad\,\,\,\text{if}~ \tau\leq 2, \\
	\left|u-u_{\lambda}^{*}\right|~\quad\quad\text{if}~ 2<\tau.
	\end{array}\right.
	\end{aligned}
\end{equation*}
Similarly on $\{v^*_\lambda\leq u\leq u^*_\lambda\}$, we have for some $c_{15}>0$
\begin{equation*}
	\left|\left|u\right|^{\tau}-\left(u_{\lambda}^{*}\right)^{\tau} \right|\leq c_{15}\left\{\begin{array}{l}
		\left| u-u^*_\lambda \right|^{\tau} \quad\,\,\,\text{if}~ \tau\leq 2,\\
		\left| u-u^*_\lambda \right|\quad\,\,\,\,\,\text{if}~ 2<\tau.
		\end{array}\right.
\end{equation*}

$\bullet$ $F(z, \cdot)$ is a $L^\infty$-locally Lipschitz integrand.\\

Using these observations in  \eqref{inequa414}, we see that for some $c_{16}>0$ and $\|u-u^*_\lambda\|\leq1$ is small
\begin{align}\label{inequa415}
\int_{\Omega}\left|G_\lambda(z, u) - G_\lambda\left(z,u_\lambda^*\right)\right| d z \leq c_{16}\left\|u-u_\lambda^*\right\|^\tau.
\end{align}

Also we have
	\begin{align}\label{inequa416}
	& \int_{\Omega}\left|G_\lambda^{+}\left(z, u_\lambda^*\right)-G_\lambda^{+}(z, u)\right| d z \notag\\
	\leq & \int_{\left\{u<v_\lambda^*\right\}}\left|G_\lambda\left(z, u_\lambda^*\right)\right| d z \notag\\
	+ & \int_{\left\{v_\lambda^* \leq u \leq u_\lambda^*\right\}}\left|\frac{\lambda}{\tau}\left(\left(u_\lambda^*\right)^\tau-\left(u^+\right)^\tau\right)+\left(F\left(z, u_\lambda^*\right)-F\left(z, u^{+}\right)\right)\right| d z \notag\\
	+ & \int_{\left\{u_\lambda^*<u\right\}} \Big\vert \lambda\left(u_\lambda^*\right)^\tau+f\left(z, u_\lambda^*\right) u_\lambda^*-\frac{\lambda}{\tau}(u^\tau-\left(u_\lambda^*\right)^\tau-(F(z, u)-F\left(z, u_\lambda^*))\right)\Big\vert d z
	\end{align}
For $\delta'\in(0,1)$, let $\bar{B}_{\delta'} (u_\lambda^*)=\Big\{u\in W_0^{1,\eta}(\Omega):\|u-u_\lambda^*\|\leq\delta'\Big\}$. Recall that $v_{\lambda}^*\prec0\prec u_{\lambda}^*$, we have
\begin{equation}\label{inequa417}
\int_{\{u<v_{\lambda}^*\}}|G_\lambda\left(z, u_\lambda^*\right)| d z \rightarrow 0, \,\,\text{and}\,\, \int_{\left\{u_\lambda^*<u\right\}}\left(\lambda\left(u_\lambda^*\right)^\tau+f\left(z, u_\lambda^*\right) u_\lambda^*\right) d z \rightarrow 0\text { as } \delta^{\prime} \rightarrow 0^{+}.
\end{equation}

Also as the above, we have for some $c_{17}>0$, all $u \in \bar{B}_{\delta^{\prime}}\left(u_\lambda^*\right)$ with $\delta^{\prime} \in(0,1)$ small
\begin{equation}\label{inequa418}
 \int_{\left\{v_\lambda^* \leq u \leq u_\lambda^*\right\}} \left|\frac{\lambda}{\tau}\left(\left(u_\lambda^*\right)^\tau-\left(u^+\right)^\tau\right)+\left(F\left(z, u_\lambda^*\right)-F(z, u^+)\right) \right| d z \leq c_{17}\left\|u-v_\lambda^*\right\|^\tau.
\end{equation}

Similarly we have for some $c_{18}>0$, all $u\in\bar{B}_{\delta^{\prime}}\left(u_\lambda^*\right)$ with $\delta^{\prime} \in(0,1)$ small,
\begin{equation}\label{inequa419}
	\int_{\left\{ u_\lambda^*<u\right\}} \left|\frac{\lambda}{\tau}\left(u^\tau-\left(u_\lambda^*\right)^\tau\right)+\left(F(z, u)-F\left(z, u_\lambda^*\right)\right) \right| dz  \leq c_{18}\|u-u_\lambda^*\|^\tau.
\end{equation}

We return to \eqref{inequa416} and use \eqref{inequa417}, \eqref{inequa418}, \eqref{inequa419}, we can obtain for some $c_{19}>0$ and for $\delta^{\prime} \in(0,1)$ small
\begin{equation}\label{inequa420}
	\int_{\Omega}\left|G_\lambda^{+}\left(z, u_\lambda^*\right)-G_\lambda^{+}(z, u)\right| d z \leq O\left(\delta^{\prime}\right)+c_{19}\left\|u-u_\lambda^*\right\|^\tau.
\end{equation}

From \eqref{inequa415} and \eqref{inequa420} it follows that
\begin{align*}
	\left|\beta_\lambda(u)-\beta_\lambda^{+}(u)\right| \leqslant c_{20} \|u-u_\lambda^* \|^\tau+O\left(\delta^{\prime}\right),
\end{align*}
for some $c_{20}>0$.

Therefore given $\varepsilon>0$, we can find $\delta_0^{\prime} \in(0,1)$ small such that
\begin{align*}
	\left|\beta_\lambda(u)-\beta_{\lambda}^{+}(u)\right| \leqslant c_{20} \delta^{\prime}+\frac{\varepsilon}{4}\quad\text{for all} ~\delta^{\prime} \in(0, \delta_0^{\prime}].
\end{align*}

Hence if $\delta_0^{\prime} \in\left(0, \frac{\varepsilon}{4 c_{20}}\right)$, we have
\begin{equation}\label{inequa421}
	\left|\beta_\lambda(u)-\beta_\lambda^{+}(u)\right| \leqslant \frac{\varepsilon}{2} \quad\text{for all} ~u \in \bar{B}_{\delta^{\prime}}\left(u_\lambda^*\right), \,\,\delta^{\prime}\in(0, \delta_0^{\prime}] .
\end{equation}

Next we estimate the difference of the two derivatives for $u,h\in W_0^{1, \eta}(\Omega)$ as follows,
\begin{align}\label{inequa422}
&\left|\langle \beta'_\lambda(u)-(\beta^+_\lambda)'(u),h\rangle \right|\notag\\
\leq&\int_{\Omega} \left|g_\lambda (z, u)-g_\lambda^{+}(z, u)\right||h| dz\notag\\
\leq&\int_{\Omega} \left|g_\lambda (z, u)-g_\lambda(z, u_\lambda^*)\right||h| dz+\int_{\Omega} \left|g_\lambda^+ (z, u_\lambda^*)-g_\lambda^{+}(z, u)\right||h| dz.
\end{align}

We know that continuous convex (and concave) functions are locally Lipschitz (see, for example,  \cite[Corollary 5.1.23]{PW}). Therefore both $g_\lambda(z,\cdot)$ and $g_\lambda^+(z,\cdot)$ are $L^\infty$-locally Lipschitz integrands (see \eqref{inequa41} and \eqref{inequa42}). So, from \eqref{inequa422} we have that for some $c_{21}>0$,
\begin{align}\label{inequa423}
	\left|\langle \beta'_\lambda(u)-(\beta^+_\lambda)'(u), h\rangle \right|\leq c_{21}\int_{\Omega}|u-u^*_\lambda|\vert h\vert dz.
\end{align}

From Proposition \ref{embeddings}, we know that $W_0^{1,\eta}(\Omega)\hookrightarrow L^{q^*}(\Omega)$ continuously. By hypotheses $(H_0)$, $2\leq q^*$, hence $(q^*)'=\frac{q^*}{q^*-1}\leq2$ and so $W_0^{1,\eta}(\Omega)\hookrightarrow L^{(q^*)'}(\Omega)$ continuously. In \eqref{inequa423}, we can use H\"older inequality and the continuous embedding $W_0^{1,\eta}(\Omega)\hookrightarrow L^{q^*}(\Omega)$, $W_0^{1,\eta}(\Omega)\hookrightarrow L^{(q^*)'}(\Omega)$ to obtain for some $c_{22}>0$
\begin{align}\label{inequa424}
& \left|\langle  \beta_\lambda^{\prime}(u)-\left(\beta_\lambda^{+}\right)^{\prime}(u), h\rangle \right| \leq c_{22}\left\|u- u_\lambda^*\right\|\|h\|\notag\\
\Rightarrow & \left\|\beta_\lambda^{\prime}(u)-\left(\beta_\lambda^{+}\right)^{\prime}(u)\right\|_* \leq c_{22}\left\|u - u_\lambda^*\right\|.
\end{align}
From \eqref{inequa421} and \eqref{inequa424} we see that given $\varepsilon>0$, we can find $\delta'>0$ such that
\begin{align*}
	\left\|\beta_\lambda-\beta_{\lambda}^+\right\|_{C^1(\bar{B}_{\delta '}(u^*_\lambda))}\leq\varepsilon.
\end{align*}
As before, the $C^{1}$-continuity property of the critical groups implies that
\begin{align*}
	C_k(\beta_\lambda,u^*_\lambda)=C_k(\beta^+_\lambda, u^*_\lambda)\quad \text{for all} ~k\in\mathbb{N}_0.
\end{align*}

In a similar fashion, we show that
\begin{align*}
	C_k(\beta_\lambda,v^*_\lambda)=C_k(\beta^-_\lambda, v^*_\lambda)\quad \text{for all} ~k\in\mathbb{N}_0.
\end{align*}

\end{proof}

Now we are ready to produce nodal solutions which vanish asymptotically as $\lambda\to0^+$.
\begin{proposition}\label{nodal}
If hypotheses $(H_0), (H_1)$ hold and $\lambda\in(0,\lambda_*)$, then problem $(P_{\lambda})$ has a nodal solution $y_\lambda\in [v^*_\lambda, u^*_\lambda]$
and
$$y_\lambda\to 0 \,\,\text{in}\,\, W_0^{1,\eta}(\Omega)\cap L^\infty(\Omega)\,\,\text{as}\,\,\lambda\to 0^+.$$
\end{proposition}

\begin{proof}
	From \eqref{inequa41} and \eqref{inequa42} it is clear that $\beta^+_\lambda(\cdot)$ is coercive. Also using Proposition \ref{embeddings}, we see that $\beta^+_\lambda(\cdot)$ is sequentially weakly lower semi-continuous. So, we can find $\widetilde{u}^*_\lambda\in W_0^{1,\eta}(\Omega) $ such that
\begin{align}\label{inequa425}
\beta^+_\lambda(\widetilde{u}^*_\lambda)=\inf \Big\{ \beta^+_\lambda(u): u\in W_0^{1,\eta}(\Omega)\Big\}.
\end{align}

	If $u\in C^1_0(\overline{\Omega} )$ with $u(z)>0$ for all $z\in\Omega$, then for $t\in(0,1)$ small, we have $0\leq tu\leq\delta$ with $\delta>0$ as in hypotheses $(H_1)(2)$. We have
	\begin{align*}
		0\leq F(z, tu(z)) \quad\text{for a.e.}~ z\in\Omega.
	\end{align*}
It follows from \eqref{inequa41}, \eqref{inequa42}, $t\in(0,1)$ and $q<p$ that
\begin{align*}
	\beta^+_\lambda(tu)\leq\frac{t^q}{q}\rho_\eta(\nabla u)-\frac{\lambda t^\tau}{\tau}\|u\|^\tau_\tau
\end{align*}

Since $\tau<q$, choosing $t\in(0,1)$ even smaller if necessary, we infer that
\begin{align*}
	&\beta^+_\lambda(tu)\leq0,\\
	\Rightarrow &\beta^+_\lambda(\widetilde{u}^*_\lambda)<0=\beta^+_\lambda(0),\\
	\Rightarrow &\widetilde{u}^*_\lambda\neq 0.
\end{align*}
From \eqref{inequa425} we have that $\widetilde{u}^*_\lambda\in K_{\beta^+_\lambda}$. Therefore \eqref{inequa43} implies that $\widetilde{u}^*_\lambda=u^*_\lambda$. So, from Proposition \ref{pro10}, we can write that
\begin{align}\label{inequa426}
	&C_k(\beta^+_\lambda,u^*_\lambda)=\delta_{k,0}\mathbb{Z} \quad\text{for all}~k\in\mathbb{N}_0,\notag\\
	\Rightarrow &C_k(\beta_\lambda,u^*_\lambda)=\delta_{k,0}\mathbb{Z} \quad\text{for all}~k\in\mathbb{N}_0.
\end{align}
Similarly, using this time $\beta_\lambda^-$, we obtain
\begin{align}\label{inequa427}
	C_k(\beta^+_\lambda,v^*_\lambda)=\delta_{k,0}\mathbb{Z} \quad\text{for all}~k\in\mathbb{N}_0.
\end{align}

From \eqref{inequa41} we see that $\beta_\lambda(\cdot)$ is coercive. Then \cite[Proposition 6.2.24]{PRR} implies that
\begin{align}\label{inequa428}
	C_k(\beta^+_\lambda,\infty)=\delta_{k,0}\mathbb{Z} \quad\text{for all}~k\in\mathbb{N}_0.
\end{align}

Suppose that $K_{\beta_\lambda}=\{0, u^*_\lambda,v^*_\lambda\}$. From Proposition \ref{crigro}, \eqref{inequa426}, \eqref{inequa427}, \eqref{inequa428} and the Morse relation with $t=-1$ (see \eqref{001}), we have
$$2(-1)^0=(-1)^0,$$
which is a contradiction. Therefore there exists $y_\lambda\in W^{1,\eta}_0(\Omega)$ such that
\begin{align*}
y_\lambda\in K_{\beta_\lambda}\subseteq [v^*_\lambda, u^*_\lambda], y_\lambda\notin \{0, u^*_\lambda, v^*_\lambda\},
\end{align*}
thus, we know that $y_\lambda$ is a nodal solution of  problem $(P_{\lambda})$ for  $\lambda\in (0,\lambda_*)$.

On account of Proposition \ref{vanishing}, we have
\begin{align*}
	& u_\lambda^*, v_\lambda^* \rightarrow 0\,\, \text {in}\,\, W_0^{1, \eta}(\Omega) \cap L^{\infty}(\Omega) \,\,\text { as }\,\, \lambda \rightarrow 0^{+}, \\
	& \Rightarrow y_\lambda \rightarrow 0 \,\,\text {in}\,\, W_0^{1, \eta}(\Omega) \cap L^{\infty}(\Omega) \,\,\text { as } \,\,\lambda \rightarrow 0^{+} .
	\end{align*}
\end{proof}

Summarizing our findings for problem $(P_{\lambda})$, we can give the proof of Theorem \ref{mainth}. Note that we provide sign information for all solutions, the solutions are ordered and we provide their asymptotic behavior as $\lambda\to0^+$.

\section{Proof of Theorem \ref{mainth}}\label{Sec5}

From Proposition \ref{vanishing}, \ref{sign} and \ref{nodal}, we know that for all $\lambda>0$ small, problem $(P_{\lambda})$ has at least three solutions $v_\lambda^*, y_\lambda, u_\lambda^* \in W_0^{1,\eta}(\Omega) \cap L^{\infty}(\Omega)$ such that
$$
\begin{gathered}
v_\lambda^* \leqslant y_\lambda \leqslant u_\lambda^*, \\
v_\lambda^* \prec 0, ~y_\lambda ~\text{is nodal},~ 0\prec u_\lambda^*, \\
\text { and } v_\lambda^*, y_\lambda, u_\lambda^*\rightarrow 0 \,\,\text {in} \,\,W_0^{1,\eta}(\Omega) \cap L^{\infty}(\Omega)\,\,\text {as} \,\,\lambda \rightarrow 0^{+}.
\end{gathered}
$$
The proof is completed.

\section*{Acknowledgements}
C. Ji was partially supported by National Natural Science Foundation of China (No. 12171152).

\noindent
\noindent \textsc{Chao Ji} \\
Department of Mathematics\\
East China University of Science and Technology \\
Shanghai 200237,  China \\
\texttt{jichao@ecust.edu.cn}\\
\noindent and \\
\textsc{Nikolaos S. Papageorgiou} \\
Department of Mathematics \\
National Technical University, Zograrou Compus \\
Athens 15780, Greece \\
\texttt{npapg@math.ntua.gr} \\




\begin{thebibliography}{90}		
\footnotesize
\bibitem{CBW}\'{A}. Crespo-Blanco, P. Winkert, Nehari manifold approach for superlinear double phase problems with variable exponents,  arXiv:2211.09189v1, 2022.

\bibitem{CGHW}\'{A}. Crespo-Blanco, L. Gasi\'{n}ski, P. Harjulehto, P. Winkert, A new class of double phase variable exponent problems: Existence and uniqueness, J. Differential Equations 323 (2022), 182-228.





\bibitem{FP} M.E. Filippakis, N.S. Papageorgiou, Multiple constant sign and nodal solutions for nonlinear elliptic equations with the $p$-Laplacian, J. Differential Equations 245 (2008),  1883-1922.


\bibitem{GP} L. Gasi\'{n}ski, N.S. Papageorgiou, Exercises in analysis. Part 2. Nonlinear Analysis,\emph{Springer}, Cham, 2016.


\bibitem{GW} L. Gasi\'{n}ski, P. Winkert, Constant sign solutions for double phase problems with superlinear nonlinearity, Nonlinear Anal. 195 (2020), 111739, 9 pp.

\bibitem{GV} M. Guedda, L. V\'{e}ron,  Quasilinear elliptic equations involving critical Sobolev exponents,  Nonlinear Anal. 13 (1989),  879¨C902.




\bibitem{HH} P. Harjulehto, P. H\"asto, Orlicz Spaces and Generalized Orlicz Spaces,  Lecture Notes in Mathematics, 2236. Springer, Cham, 2019.




\bibitem{HP} S.C. Hu, N.S. Papageorgiou, Research topics in analysis, Vol. I. Grounding theory,  Birkh\"{a}user Advanced Texts: Basler Lehrb\"{u}cher, Birkh\"{a}user/Springer, Cham, 2022.


\bibitem{LD} W.L. Liu, G.W. Dai, Existence and multiplicity results for double phase problem, J. Differential Equations 265 (2018), 4311-4334.



\bibitem{LP} Z.H. Liu, N.S. Papageorgiou,  On an anisotropic eigenvalue problem, Results Math. 78 (2023),  Paper No. 178.



\bibitem{Ma} P. Marcellini,   Regularity of minimizers of integrals of the calculus of variations with nonstandard growth conditions, Arch. Rational Mech. Anal. 105 (1989),  267-284




\bibitem{Ma1} P. Marcellini, Regularity and existence of solutions of elliptic equations with $p,q$-growth conditions, J. Differential Equations 90 (1991), 1-30.


\bibitem{Ma2} P. Marcellini, Growth conditions and regularity for weak solutions to nonlinear elliptic pdes,  J. Math. Anal. Appl. 501 (2021), 124408, 32 pp.


\bibitem{MW} G.  Marino, P. Winkert,   Moser iteration applied to elliptic equations with critical growth on the boundary,  Nonlinear Anal. 180 (2019), 154-169.




\bibitem{MR} G. Mingione, V. R\v{a}dulescu, Recent developments in problems with nonstandard growth and nonuniform ellipticity, J. Math. Anal. Appl. 501 (2021),  125197, 41 pp.




\bibitem{MP} D. Mugnai, N.S. Papageorgiou, Resonant nonlinear Neumann problems with indefinite weight,  Ann. Sc. Norm. Super. Pisa Cl. Sci.  11 (2012), 729¨C788.



\bibitem{PRR} N.S. Papageorgiou, V.D. R\v{a}dulescu, D.D. Repov\v{s}, Nonlinear analysis¡ªtheory and methods, Springer Monographs in Mathematics, Springer, Cham, 2019.


\bibitem{PVV} N.S.  Papageorgiou, C. Vetro, F. Vetro, Multiple solutions for parametric double phase Dirichlet problems,
Commun. Contemp. Math. 23 (2021),  2050006, 18 pp.




\bibitem{PW} N.S. Papageorgiou,  P. Winkert,  Applied Nonlinear Functional Analysis. An introduction,  De Gruyter Graduate. De Gruyter, Berlin, 2018.




\bibitem{PZ} N.S. Papageorgiou, C. Zhang, Multiple ground-state solutions with sign information for double-phase Robin problems,
Israel J. Math. 253 (2023), 419-443.


\bibitem{Zh1} V.V. Zhikov, Averaging of functionals of the calculus of variations and elasticity theory, Math. USSR Izv., 29 (1987), 33-66




\bibitem{Zh2} V.V. Zhikov,  On Lavrentiev's phenomenon,  Russian J. Math. Phys. 3 (1995), 249-269.



\end{thebibliography}
\end{document}